\documentclass[a4paper]{article}

\usepackage{amsmath}
\usepackage{amsfonts}
\usepackage{amssymb}
\usepackage{amsthm}
\usepackage{color}
\usepackage{graphicx}
\usepackage{float,wrapfig}
\usepackage{url}
\usepackage[latin1]{inputenc}

\title{{\large \textbf{Amenable actions of amalgamated free products}}}
\author{Soyoung Moon\footnote{This work is supported by Swiss NSF grant $20\_118014/1$.}}
\date{\today}

\newcommand{\supp}{\textrm{supp}}
\newcommand{\Fix}{\textrm{Fix}}

\theoremstyle{plain}
\newtheorem{thm}{Theorem}
\newtheorem{prop}[thm]{Proposition}
\newtheorem{cor}[thm]{Corollary}
\newtheorem{lem}[thm]{Lemma}

\theoremstyle{definition}
\newtheorem{defn}{Definition}[section]

\newtheorem{ex}{Example}[section]
\begin{document}
\maketitle

\begin{abstract}
We prove that the amalgamated free product of two free groups of rank two over a common cyclic subgroup, admits an amenable, faithful, transitive action on an infinite countable set. We also show that any finite index subgroup admits such an action, which applies for example to surface groups and fundamental groups of surface bundles over $\mathbb{S}^1$.
\end{abstract}
\pagestyle{headings}

\section{Introduction}

An action of a group $G$ on a set $X$ is \textit{amenable} if there exists a $G$-invariant mean on $X$, i.e. a map $\mu : 2^X=\mathcal{P}(X)
\rightarrow [0,1]$ such that $\mu(X)=1 $, $\mu (A \cup B)= \mu(A)+\mu(B)$ for every pair of disjoint subsets $A$, $B$ of $X$, and
$\mu(gA)=\mu(A)$, $\forall g \in G$, $\forall A \subseteq X$.

\smallskip

The study of amenability goes back to von Neumann \cite{vN} and has spanned over the 20th century in various fields of mathematics, such as geometric group theory, harmonic analysis, graph theory, operator algebra, etc. F. P. Greenleaf asked in \cite{Green} whether the presence of a $G$-invariant mean on a set on which $G$ acts faithfully implies that the group $G$ is amenable (i.e. if the action on itself by left multiplication is amenable), and the first counter example was given in \cite{van}, where E. K. van Douwen constructed an interesting amenable action of the non-abelian free group.

The above definition is due to Greenleaf \cite{Green}. We should mention that Zimmer \cite{Zim} has also introduced a notion of amenability for a group action that is different from ours; an action by homeomorphisms of a countable discrete group $G$ on a compact Hausdorff space $X$ is \textit{(topologically) Zimmer amenable} if there exists a sequence of continuous maps $m^n:X\rightarrow Proba(G)$ such that $\lim_{n\rightarrow \infty} \textrm{sup}_{x\in X} \|gm^n_x - m^n_{gx}\|_1 =0$, for all $g\in G$ (cf. \cite{Ozawa}, \cite{HgRoe}, \cite{AnRe}). With this definition, a group is amenable if and only if the action on an one-point space is Zimmer amenable, while such an action is always Greenleaf amenable. On the other hand, the action of $G$ on itself by left multiplication is always Zimmer amenable (by taking $m^n : G \rightarrow Proba(G)$ defined by $m^n_g=\delta_g$). More generally, the action of $G$ on a homogenous space $G/H$ is Zimmer amenable if and only if the subgroup $H$ is amenable. From now on, we will use the term of an amenable action as mean of Greenleaf amenable action.

\smallskip

For the study of amenable actions of a group $G$, we should require some restrictions on the $G$-action in order to avoid trivial cases. One should assume that the action is faithful, otherwise one would take immediately a free group $\mathbb{F}_n$, $n\geq 2$, and any non-trivial normal subgroup $N \triangleleft \mathbb{F}_n$ such that the quotient group $\mathbb{F}_n / N$ is amenable (e.g. $N=\mathbb{F}_n'$ the commutator subgroup), so that the natural action of $\mathbb{F}_n$ on $\mathbb{F}_n / N$ is amenable but not faithful. In addition, one should require that $G$ acts transitively, otherwise one could take any group $G$ and $X=G\sqcup Y$ where $G$ acts on $Y$ amenably, so that the $G$-action on $X$ is faithful and amenable (since there is a $G$-equivariant map from $Y$ into $X$). In this direction, Y. Glasner and N. Monod \cite{GlMo} proposed to study the class $\mathcal{A}$ of all countable groups which admit a faithful, transitive and amenable action. The class $\mathcal{A}$ is closed under direct products and free products, and a group is in $\mathcal{A}$ if it has a co-amenable subgroup which is in $\mathcal{A}$ (Proposition 1.7 in \cite{GlMo}). On the other hand, in general the class is neither closed under passing to subgroups (the case of finite index subgroups is open), nor closed under semidirect products. As an example for semidirect product, one may take the group $SL_2(\mathbb{Z}) \ltimes \mathbb{Z}^2$; while $SL_2(\mathbb{Z})$ is in $\mathcal{A}$ since it contains a free group of finite index, the pair $(SL_2(\mathbb{Z}) \ltimes \mathbb{Z}^2, \mathbb{Z}^2)$ has the relative property (T) (cf \cite{BHV}), so that the group $SL_2(\mathbb{Z}) \ltimes \mathbb{Z}^2$ is not in $\mathcal{A}$ (Lemma 4.3 in \cite{GlMo}). Besides, this group is another example which shows that the class $\mathcal{A}$ is not closed under amalgamated free products; one may see the group $SL_2(\mathbb{Z}) \ltimes \mathbb{Z}^2$ as the amalgamated free product $G\ast_A H$ of $G=\mathbb{Z}/4\mathbb{Z} \ltimes \mathbb{Z}^2$ and $H=\mathbb{Z}/6\mathbb{Z} \ltimes \mathbb{Z}^2$ along $A=\mathbb{Z}/2\mathbb{Z} \ltimes \mathbb{Z}^2$ and notice that the three groups $G$, $H$ and $A$ are in $\mathcal{A}$ since they are amenable.

In particular, Y. Glasner and N. Monod showed that the free product of any two countable groups is in $\mathcal{A}$ unless one factor has the fixed point property and the other has the virtual fixed point property\footnote{A group $G$ has the \textit{fixed point property} if any amenable $G$-action has a fixed point, and $G$ has the \textit{virtual fixed point property} if it has a finite index subgroup having the fixed point property.}; for this, they used an argument of genericity in Baire's sense (Theorem 3.3 in \cite{GlMo}). Let us mention that another construction of amenable action of a non-abelian free group is obtained by R. Grigorchuk and V. Nekrashevych in \cite{GrNe}.

\smallskip

The main result of this paper is, motivated by this method of genericity, to give another example of non-amenable group which is in $\mathcal{A}$ (see Theorem \ref{amalgam} and Theorem \ref{amalgam_any_c}):

\bigskip

\noindent \textbf{Theorem.} The amalgams $\mathbb{F}_2 \ast_\mathbb{Z} \mathbb{F}_2$ belong to $\mathcal{A}$, where $\mathbb{Z}$ embeds in each factor as subgroup generated by some common word on the generating sets.

\bigskip
Such amalgams are known as doubles of $\mathbb{F}_2$.
The key point of the proof is to fix a transitive permutation $\beta$ and take a generic element $\alpha$ (i.e. an element in the intersection of countably many generic sets) in order to construct $\mathbb{F}_2=\langle \alpha, \beta \rangle$ in a way that the amalgamated free product of two copies of $\mathbb{F}_2$ along a cyclic group has the desired properties. Therefore, the difficulty of the proof resides in the choice of the generic sets because they can be very ``nasty" (see Proposition \ref{F_2-fidel1}).

As we mentioned before, in general it is not known whether the class $\mathcal{A}$ is closed under passing to finite index subgroups or not. But it is true for our case (see Theorem \ref{amalgam-index-in-A}):

\bigskip

\noindent \textbf{Theorem.} For any finite index subgroup $H$ of $\mathbb{F}_2 \ast_\mathbb{Z} \mathbb{F}_2$ as above, $H$ belongs to $\mathcal{A}$.

\bigskip

A surface group $\Gamma_g$ is the fundamental group of a closed oriented surface of genus $g\geq 2$. The group $\Gamma_2$ can be viewed as an amalgamated free product of two copies of $\mathbb{F}_2$ along the subgroup generated by the commutator, i.e. $\Gamma_2= \langle a_1,b_1\rangle \ast_{\langle c \rangle} \langle a_2,b_2\rangle$ where $c=[a_1,b_1]=[a_2,b_2]$. For $g\geq 3$, $\Gamma_g$ injects into $\Gamma_2$ as a finite index subgroup. Therefore, by applying our results, we have the following theorem (see Theorem \ref{Sigma_g}):

\bigskip

\noindent \textbf{Theorem.} The surface groups $\Gamma_g$ belong to $\mathcal{A}$, $\forall g\geq 2$.

\bigskip

As a corollary, we obtain that the fundamental group of a 3-manifold which virtually fibers over the circle is in $\mathcal{A}$. Indeed, let $M$ be a 3-manifold which fibers over the circle. Then there is a short exact sequence:
$$
0 \rightarrow \Gamma_g \rightarrow \pi_1(M) \rightarrow \mathbb{Z} \rightarrow 0,
$$
so that the subgroup $\Gamma_g$ is co-amenable in $\pi_1(M)$. Moreover, if $M$ is a 3-manifold which virtually fibers over the circle, then it contains a finite index subgroup which is in $\mathcal{A}$, so that $\pi_1(M)$ is also in $\mathcal{A}$. Some examples of the fundamental group of such manifolds are given in \cite{agol}, which includes the Bianchi groups PSL$(2,\mathcal{O}_d)$, where $\mathcal{O}_d$ is the ring of integers of the imaginary quadratic field $\mathbb{Q}(\sqrt{-d})$ with $d$ a positive integer.

\bigskip
\bigskip

\noindent \textbf{Acknowledgement.} I would like to thank Nicolas Monod for suggesting the question and for helpful discussions, Alain Valette for his constant help and encouragement, and the referee for useful comments on the first version of this paper.

\section{Baire spaces}

For the importance of the idea of generic choice, we briefly discuss Baire spaces in this chapter.
\begin{defn}
A topological space $X$ is a \textit{Baire space} if every intersection of countably many dense open subsets is dense in $X$.
\end{defn}
Equivalently, $X$ is a Baire space if every union of countably many closed subsets with empty interior has empty interior.

\begin{defn}
A \textit{Polish space} is a separable completely metrizable topological space, i.e. it is a space homeomorphic to a complete space that has a countable dense subset.
\end{defn}

Observe that any closed subspace of a Polish space is Polish.

\bigskip

Let $X$ be an infinite countable set. Equipped with the discrete topology, $X$ is a complete topological space. Let us denote by $X^X$ the set of all self-maps of $X$ and endow it with the topology of pointwise convergence (i.e. $\alpha_n$ converges to $\alpha$ if for all finite subset $F$ of $X$, there exists $n_0$ such that $\alpha_n|_F=\alpha |_F$, for all $n \geq n_0$). This is the product of the topologies of $X$. Hence $X^X$ is complete being a product of complete spaces, and it is separable and metrizable since it is a countable product of separable, metrizable spaces. So $X^X$ is a Polish space and by Baire's theorem, it is a Baire space.

Let us denote by $Sym(X)\subset X^X$ the group of permutations of $X$. Equipped with the induced topology of $X^X$, $Sym(X)$ is a topological group. Indeed, let $\{\alpha_n\}_{n\geq 1}$ be a sequence converging to $\alpha$ in $Sym(X)$. Let $F\subset X$ be a finite subset of $X$. There exists $n_0$ such that $\alpha_n |_{F\cup \alpha^{-1}F}=\alpha |_{F\cup \alpha^{-1}F}$, $\forall n\geq n_0$. Then for all $x\in F$, we have $\alpha_n(\alpha^{-1}(x))=\alpha(\alpha^{-1}(x))=x$, so $\alpha_n^{-1}(x)=\alpha^{-1}(x)$, $\forall n\geq n_0$. Therefore $\alpha_n^{-1}$ converges to $\alpha^{-1}$, so that the application $\alpha \mapsto \alpha^{-1}$ is continuous. Moreover, let $\{\beta_m\}_{m\geq 1}$ be a sequence converging to $\beta$ in $Sym(X)$. Let $F \subset X$ be a finite subset of $X$. There exists $n_1$ such that $\alpha_n |_{F\cup \beta F}=\alpha |_{F\cup \beta F}$, $\forall n\geq n_1$. In addition, there exists $n_2$ such that $\beta_m |_F=\beta |_F$, for all $m\geq n_2$. Then for all $x\in F$, $\alpha_n(\beta_m(x))=\alpha_n(\beta(x))=\alpha\beta(x)$, for all $m \geq \mathrm{max}\{n_1,n_2\}$. Therefore $\alpha_n \beta_m$ converges to $\alpha\beta$, so that the application $(\alpha,\beta) \mapsto \alpha\beta$ is continuous.

Consequently, the injection $i:Sym(X) \rightarrow X^X \times X^X$; $\alpha \mapsto (\alpha, \alpha^{-1})$ is a homeomorphism onto its image which is closed. Thus $Sym(X)$ is a Polish space, in particular it is a Baire space.

\begin{defn}
A subset $Y\subset Sym(X)$ is called
\begin{enumerate}
\item[$\cdot$] \textit{meagre} or \textit{first category} if it is a union of countably many closed subsets with empty interior;
\item[$\cdot$] \textit{generic} or \textit{dense $G_{\delta}$} if its complement $Sym(X)\setminus Y$ is meagre, i.e. it is an intersection of countably many dense open subsets.
\end{enumerate}
\end{defn}

By definition of the topology on $Sym(X)$, a subset $Y\subset Sym(X)$ has empty interior if for all $\alpha' \in Y$ and for all finite subset $F\subset X$, there exists $\alpha\in Sym(X)\setminus Y$ such that $\alpha'|_F=\alpha|_F$.

\section{Construction of $\mathbb{F}_2$}
Let $X$ be an infinite countable set. Let $\beta$ be a simply transitive permutation of $X$. Let $c=c(\alpha,\beta)$ be a weakly cyclically reduced word (i.e. if $c=g_m\cdots g_1$, then $g_m \neq g_1^{-1}$) on the alphabet $\{\alpha^{\pm 1}$, $\beta^{\pm 1}\}$ such that $c\notin \langle\beta\rangle$.

\begin{prop}\label{F_2-fidel1}
The set
\begin{eqnarray}
\mathcal{U}_1=\{ \alpha \in Sym(X) & | & \forall w\in \langle \alpha, \beta \rangle \setminus \langle c \rangle, \textrm{ there exist infinitely many $x\in X$} \nonumber\\
&& \textrm{such that $cx=x$, $cwx=wx$ and $wx \neq x$ } \} \nonumber
\end{eqnarray}
is generic in $Sym(X)$.
\end{prop}

\begin{prop}\label{F_2-fidel2}
The set
$$
\mathcal{U}_2 = \{\alpha\in Sym(X) \mid \forall k \in \mathbb{Z}\setminus \{0\}, \exists x\in X \textrm{ such that } c^kx\neq x\}
$$
is generic in $Sym(X)$.
\end{prop}

Note that $\mathcal{U}_2$ is the set of $\alpha$'s such that $c$ has infinite order.

\begin{defn}
Let $c=c(\alpha, \beta)$ be a weakly cyclically reduced word. Let $S(\alpha)$ be the sum of exponents of $\alpha$, and $S(\beta)$ be the sum of exponents of $\beta$. We say that $c$ is \emph{special} if $c$ is one of the following types:
\begin{enumerate}
\item[(1)] $S(\alpha)=S(\beta)=0$;
\item[(2)] $S(\alpha)$ divides $S(\beta)$.
\end{enumerate}
\end{defn}

Let $\{A_n\}_{n=1}^{\infty}$ be a pairwise disjoint F\o lner sequence for $\beta$, that is
$$ \lim_{n\rightarrow \infty}\frac{|A_n \vartriangle \beta \cdot A_n|}{|A_n|}=0. $$

\begin{prop}\label{F_2-moy}
Let $c$ be a special word. The set
\begin{eqnarray}
\mathcal{U}_3 =\{ \alpha \in Sym(X) &|& \textrm{there exists $\{A_{n_k}\}_{k=1}^{\infty}$ a subsequence of $\{A_n\}_{n=1}^{\infty}$} \nonumber \\
&& \textrm{such that $A_{n_k}\subset \mathrm{Fix}(c)$, $\forall k\geq 1$ and $\{A_{n_k}\}_{k=1}^{\infty}$ is a} \nonumber \\ 
&& \textrm{F\o lner sequence for $\alpha$}\} \nonumber
\end{eqnarray}
is generic in $Sym(X)$.
\end{prop}

\begin{prop}\label{F_2-transitive}
The set
\begin{eqnarray}
\mathcal{U}_4=\{\alpha\in Sym(X)&|&\textrm{for every finite index subgroup $H$ of $\langle \alpha,\beta\rangle$, the action}\nonumber \\
 &&\textrm{of $H$ on $X$ is transitive } \} \nonumber
\end{eqnarray}
is generic in $Sym(X)$.
\end{prop}

From the previous four propositions, one deduces immediately:
\begin{cor}\label{F_2-corollary}
Let $c$ be a special word on $\{\alpha^{\pm 1}$, $\beta^{\pm 1}\}$. Let $\alpha\in \mathcal{U}_1 \cap \mathcal{U}_2 \cap \mathcal{U}_3 \cap \mathcal{U}_4$. Then $\langle \alpha$, $\beta \rangle \simeq \mathbb{F}_2$ and
\begin{enumerate}
\item[(1)] the action of $\mathbb{F}_2$ on $X$ is transitive and faithful;
\item[(2)] for all $w\in \langle \alpha$, $\beta \rangle \setminus \langle c \rangle$, there exist infinitely many $x\in X$ such that $cx=x$, $cwx=wx$ and $wx \neq x$. In particular, there are infinitely many fixed points of $c$ in $X$;
\item[(3)] there exists a pairwise disjoint F\o lner sequence for $\langle \alpha$, $\beta \rangle$ which is fixed by $c$;
\item[(4)] for all finite index subgroup $H$ of $\langle \alpha,\beta\rangle$, the $H$-action on $X$ is transitive.
\end{enumerate}

\end{cor}

\subsection{Proofs of Propositions \ref{F_2-fidel1} and \ref{F_2-fidel2}}

Propositions \ref{F_2-fidel1} and \ref{F_2-fidel2} are sufficient conditions for faithfulness of $\mathbb{F}_2$-action with some additional ``unnatural looking'' properties that will be needed for construction of $\mathbb{F}_2 \ast_\mathbb{Z} \mathbb{F}_2$ in Chapter 4. As we resort to the graph theory for these proofs, we begin by fixing the notations on graphs that will be used in the section. The fundamental notions are based on \cite{Serre}.

\subsubsection{Graph extension}

A graph $G$ consists of the set of vertices $V(G)$ and the set of edges $E(G)$, and two applications $E(G) \rightarrow E(G)$; $e \mapsto \bar{e}$ such that $\bar{\bar{e}}=e$ and $\bar{e}\neq e$, and $E(G)\rightarrow V(G)\times V(G)$; $e \mapsto (i(e), t(e))$ such that $i(e)=t(\bar{e})$. An element $e\in E(G)$ is a \textit{directed edge} of $G$ and $\bar{e}$ is the \textit{inverse edge} of $e$. For all $e\in E(G)$, $i(e)$ is the \textit{initial vertex} of $e$ and $t(e)$ is the \textit{terminal vertex} of $e$.

Let $S$ be a set. A \textit{labeling} of a graph $G=(V(G), E(G))$ on the set $S^{\pm 1}=S\cup S^{-1}$ is an application
$$
l:E(G) \rightarrow S^{\pm 1}; e \mapsto l(e)
$$
such that $l(\bar{e})=l(e)^{-1}$. A \textit{labeled graph} $G=(V(G), E(G),S,l)$ is a graph with a labeling $l$ on the set $S^{\pm 1}$. A labeled graph is \textit{well-labeled} if for any edges $e$, $e' \in E(G)$, $\big[ i(e)=i(e')$ and $l(e)=l(e')\big]$ implies that $e=e'$. If a group $\Gamma=\langle S\rangle$ acts on $X$, a labeled graph with set of vertices $X$ and set of edges $S^{\pm 1}$ is well-labeled if and only if it is a Schreier graph.

\bigskip

A word $w=w_m\cdots w_1$ on $\{ \alpha^{\pm 1}, \beta^{\pm 1} \}$ is called \textit{reduced} if $w_{k+1}\neq w_k^{-1}$, $\forall 1\leq k\leq m-1$. A word $w=w_m\cdots w_1$ on $\{ \alpha^{\pm 1}, \beta^{\pm 1} \}$ is called \textit{weakly cyclically reduced} if $w$ is reduced and $w_m\neq w_1^{-1}$; this definition allows that $w_m$ and $w_1$ to be equal. We denote by $|w|$ the word length of $w$. Given a reduced word, we shall define two finite graphs labeled on $\{ \alpha^{\pm 1}, \beta^{\pm 1} \}$ as follows:

\begin{defn}\label{def-path}
Let $w=w_m\cdots w_1$ be a reduced word on $\{ \alpha^{\pm 1}, \beta^{\pm 1} \}$. The \emph{path} of $w$ is a finite labeled graph $P(w, v_0)$ consisting of $|w|+1$ vertices and $|w|$ directed edges $\{e_1$, $\dots$, $e_m\}$ such that
\begin{enumerate}
\item[$\cdot$] $i(e_{j+1})=t(e_j)$, $\forall 1 \leq j \leq m-1$;
\item[$\cdot$] $v_0=i(e_1)\neq t(e_m)$;
\item[$\cdot$] $l(e_j)=w_j$, $\forall 1\leq j \leq m$.
\end{enumerate}
\end{defn}

\begin{figure}[H]
\centering
\includegraphics[width=10.5cm]{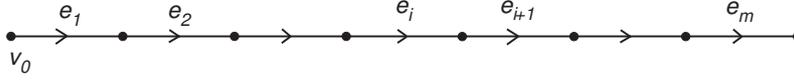}
\caption{The path of $w$}
\label{}
\end{figure}

\begin{defn}\label{def-cycle}
Let $w=w_m\cdots w_1$ be a reduced word on $\{ \alpha^{\pm 1}, \beta^{\pm 1} \}$. The \emph{cycle} of $w$ is a finite labeled graph $C(w, v_0)$ consisting of $|w|$ vertices and $|w|$ directed edges $\{e_1$, $\dots$, $e_m\}$ such that
\begin{enumerate}
\item[$\cdot$] $i(e_{j+1})=t(e_j)$, $\forall 1 \leq j \leq m-1$;
\item[$\cdot$] $v_0=i(e_1)= t(e_m)$;
\item[$\cdot$] $l(e_j)=w_j$, $\forall 1\leq j \leq m$.
\end{enumerate}
\end{defn}

\begin{figure}[H]
\centering
\includegraphics[width=3.5cm]{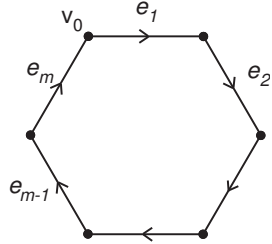}%
\caption{The cycle of $w$}
\label{path-basic}
\end{figure}

Notice that since $w$ is a reduced word, the graph $P(w, v_0)$ is well-labeled. If $w$ is weakly cyclically reduced, then $C(w, v_0)$ is also well-labeled.

Reciprocally, if $P=\{e_1$, $e_2$, $\dots$, $e_n\}$ is a well-labeled path with $i(e_1)=v_0$, labeled by $l(e_i)=g_i$, $\forall i$, then there exists a unique reduced word $w=g_n\cdots g_1$ such that $P(w,v_0) $ is $P$. If $C=\{e_1$, $e_2$, $\dots$, $e_n\}$ is a well-labeled cycle with $t(e_n)=i(e_1)=v_0$, labeled by $l(e_i)=g_i$, $\forall i$, then there exists a unique weakly cyclically reduced word $w=g_n\cdots g_1$ such that $C(w,v_0)$ is $C$.

\bigskip

Let $X$ be an infinite countable set. Let $\beta$ be a simply transitive permutation of $X$. We shall represent the $\beta$-action on $X$ as an infinite 2-regular well-labeled graph. The \textit{pre-graph} $G_0$ is a labeled graph consisting of the set of vertices $V(G_0)=X$ and the set of edges $E(G_0)$ where for all $e \in E(G_0)$, $l(e)\in \{\beta^{\pm 1 }\}$ and such that every vertex has exactly one entering edge and one leaving edge. One can imagine $G_0$ as the Cayley graph of $\mathbb{Z}$ with $1$ as a generator.

\begin{defn}
An \textit{extension} of $G_0$ is a well-labeled graph $G$ labeled by $\{\alpha^{\pm 1}, \beta^{\pm 1}\}$, containing $G_0$. We will denote it by $G_0 \subset G$.
\end{defn}

\begin{figure}[H]
\centering
\includegraphics[width=3.5cm]{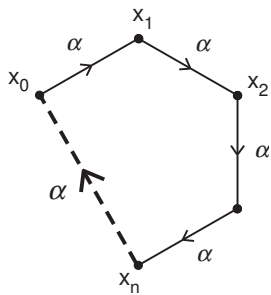}
\caption{The $\alpha$-orbit of $x_0$ that has the size $n+1$.}
\label{buckle}
\end{figure}

In order to have a transitive action with some additional properties of the $\langle \alpha, \beta \rangle$-action on $X$, we shall extend $G_0$ by adding finitely many directed edges labeled by $\alpha$ on $G_0$ where the edges labeled by $\beta$ are already prescribed. In order that the added edges represent an action on $X$, we put the edges in such a way that the extended graph is well-labeled, and moreover we put an additional edge labeled by $\alpha$ on every endpoint of the extended edges by $\alpha$; more precisely, if we have added $n$ edges labeled by $\alpha$ between $x_0$, $x_1$, $\dots$, $x_n$ successively, we put an $\alpha$-edge from $x_n$ to $x_0$ to have a cycle consisting of $n+1$ edges (see Figure \ref{buckle}). On the points where no $\alpha$-edges are involved, we put a loop labeled by $\alpha$; this means that these points are the fixed points of $\alpha$. In the end, every point has a entering edge and a leaving edge labeled by $\alpha$ (the entering edge is equal to the leaving edge if the edge is a loop), so that the graph represents an $\langle \alpha, \beta \rangle$-action on $X$, and every $\alpha$-orbit is finite.

\begin{defn}
Let $G$, $G'$ be graphs labeled by $\{\alpha^{\pm 1}$, $\beta^{\pm 1} \}$. A \textit{homomorphism} $f: G \rightarrow G'$ is a map sending vertices to vertices, edges to edges, such that
\begin{enumerate}
\item[$\cdot$] $f(i(e))=i(f(e))$ and $f(t(e))=t(f(e))$;
\item[$\cdot$] $l(e)=l(f(e))$,
\end{enumerate}
for all $e\in E(G)$.
\end{defn}
If there exists an injective homomorphism $f:G \rightarrow G'$, we say that $f$ is an \textit{embedding}, and $G$ \textit{embeds} in $G'$. If there exists a bijective homomorphism $f:G \rightarrow G'$, we say that $f$ is an \textit{isomorphism}, and $G$ is \textit{isomorphic} to $G'$.

\begin{prop}\label{prop-path}
Let $w=w_m\cdots w_1$ be a reduced word on $\{\alpha^{\pm 1}$, $\beta^{\pm 1} \}$. Let $P(w, v_0)=\{e_1,\dots, e_m\}$ be the path defined in Definition \ref{def-path}.
There exists an extension $G$ of $G_0$ such that $P(w,v_0)$ embeds in $G$, and $P(w,v_0)$ is isomorphic to its image by the corresponding embedding. In particular, the image of $P(w,v_0)$ is a path in $G$.
\end{prop}

\begin{figure}
\centering
\includegraphics[width=5.5cm]{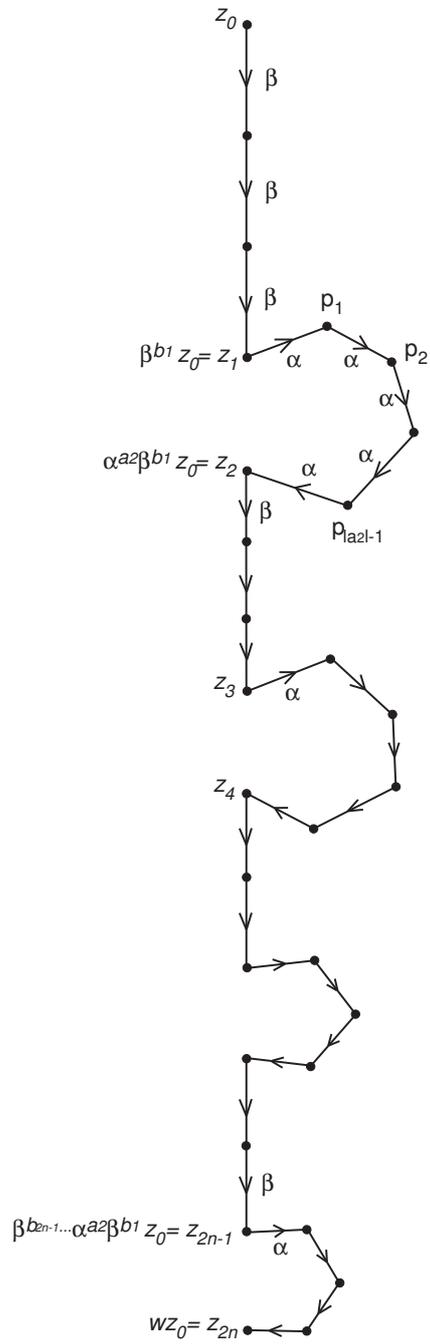}
\caption{Construction of a path in $G$}
\label{properpath}
\end{figure}

\begin{proof}
 It is enough to consider the case where $w=\alpha^{a_{2n}}\beta^{b_{2n-1}}\cdots\alpha^{a_4} \beta^{b_3} \alpha^{a_2}\beta^{b_1}$, with $m=\sum_{i=1}^n \big(|b_{2i-1}|+|a_{2i}|\big)$. Indeed, the other three cases follow from this case by taking $n$ large enough since we are treating all subwords of $w$. Let $N=\textrm{max}_j |b_j|$. For $z\in X$, denote by $B_N(z)=\{\beta^{l}z \mid -N \leq l \leq N\}$ a segment in the $\beta$-orbit of $z$.

Choose $z_0 \in X$. For all $1\leq k \leq n$, we extend $G_0$ inductively by applying the following algorithm:

\bigskip

\noindent \textbf{Algorithm (A)}
\begin{enumerate}
\item[(1)] Let $z_{2k-1}=\beta^{b_{2k-1}}z_{2k-2}$;
\item[(2)] Choose $z_{2k}\in X$ such that $B_N(z_{2k})$ is outside of the finite set of all used points;
\item[(3)] Choose $|a_{2k}|-1$ points $\{ p^{(a_{2k})}_1$, $\dots$, $p^{(a_{2k})}_{|a_{2k}|-1}\}$ outside of the finite set of all points used so far;
\item[(4)] Put the directed edges labeled by $\alpha^{sign(a_{2k})}$ from
\begin{enumerate}
\item[$\cdot$] $z_{2k-1}$ to $p^{(a_{2k})}_1$;
\item[$\cdot$] $p^{(a_{2k})}_j$ to $p^{(a_{2k})}_{j+1}$, $\forall 1 \leq j \leq |a_{2k}|-2$;
\item[$\cdot$] $p^{(a_{2k})}_{|a_{2k}|-1}$ to $z_{2k}$,
\end{enumerate}
so that we have $\alpha^{a_{2k}}z_{2k-1}=z_{2k}$.
\end{enumerate}
In the end, we have added $\sum_{i=1}^n |a_{2i}|$ new directed edges labeled by $\alpha$ (or $\alpha^{-1}$) on $G_0$ (see Figure \ref{properpath}). Let $G$ be the extended graph of $G_0$. In this construction, we have considered $|w|+1$ points $\big\{z_0$, $\beta^{sign(b_1)}z_0$, $\beta^{2sign(b_1)}z_0$, $\dots$, $\beta^{b_1}z_0 =z_1$, $\alpha^{sign(a_2)}\beta^{b_1}z_0$, $\dots$, $\alpha^{a_2}\beta^{b_1}z_0=z_2$, $\dots$, $wz_0\big\}$ in $X$, that are
$$
\textrm{$\{z_0$, $w_1z_0$, $w_2 w_1 z_0$, $\dots$, $wz_0\}$}
$$
with $l\big((w_{k-1}\cdots w_1 z_0), (w_k w_{k-1}\cdots w_1 z_0)\big)=w_k$, where $(p_1, p_2)$ symbolizes the edge $e$ with $i(e)=p_1$ and $t(e)=p_2$.

Now, we define an embedding $f:P(w,v_0)\hookrightarrow G$ by
\begin{eqnarray}
E(P(w,v_0)) &\rightarrow& E(G) \nonumber \\
e_1=(v_0, t(e_1)) &\mapsto & (z_0, w_1 z_0), \nonumber \\
e_k=(i(e_k), t(e_k)) &\mapsto & (w_{k-1}\cdots w_1 z_0, \, w_k\cdots w_1 z_0), \quad \forall 2\leq k \leq m. \nonumber
\end{eqnarray}
By construction, $P(w,v_0)$ is isomorphic to its image.
\end{proof}

\begin{prop}\label{prop-cycle}
Let $w=w_m\cdots w_1$ be a weakly cyclically reduced word on $\{\alpha^{\pm 1}$, $\beta^{\pm 1} \}$ with $w \notin \langle \beta \rangle$. Let $C(w, v_0)=\{e_1,\dots,e_m\}$ be the cycle defined in Definition \ref{def-cycle}.
There exists an extension $G$ of $G_0$ such that $C(w,v_0)$ embeds in $G$, and $C(w,v_0)$ is isomorphic to its image by the corresponding embedding. In particular, the image of $C(w,v_0)$ is a cycle in $G$.
\end{prop}

\begin{proof}
It is enough to consider the case where $w=\alpha^{a_{2n}}\beta^{b_{2n-1}}\cdots\alpha^{a_4} \beta^{b_3} \alpha^{a_2}\beta^{b_1}$, with $m=\sum_{i=1}^n \big(|b_{2i-1}|+|a_{2i}|\big)$. Let $N=\textrm{max}_j |b_j|$.

Choose $z_0 \in X$. We extend $G_0$ inductively by applying Algorithm (A) for $1\leq k \leq n-1$. Let $z_{2n-1}=\beta^{b_{2n-1}}z_{2n-2}$. Choose $|a_{2n}|-1$ points $\{ p_1$, $\dots$, $p_{|a_{2n}|-1}\}$ outside of the finite set of all points used so far. Put the directed edges labeled by $\alpha^{sign(a_{2n})}$ from
\begin{enumerate}
\item[$\cdot$] $z_{2n-1}$ to $p_1$;
\item[$\cdot$] $p_j$ to $p_{j+1}$, $\forall 1 \leq j \leq |a_{2n}|-2$;
\item[$\cdot$] $p_{|a_{2n}|-1}$ to $z_{0}$.
\end{enumerate}

We define an embedding $f:C(w,v_0)\hookrightarrow G$ by
\begin{eqnarray}
E(C(w,v_0)) &\rightarrow& E(G) \nonumber \\
e_1=(v_0, t(e_1)) &\mapsto & (z_0,\, w_1 z_0), \nonumber \\
e_k=(i(e_k), t(e_k)) &\mapsto & (w_{k-1}\cdots w_1 z_0, \, w_k\cdots w_1 z_0), \quad \forall 2\leq k \leq m-1, \nonumber \\
e_m=(i(e_m), v_0) &\mapsto & (w_{m-1}\cdots w_1 z_0,\,z_0). \nonumber
\end{eqnarray}
By construction, $C(w,v_0)$ is isomorphic to its image.

\end{proof}

\begin{cor}\label{cor-path-cycle}
Let $w$ be a reduced word. Let $F\subset G_0$ be a finite subset of $X$. There exists an extension $G$ of $G_0$ such that $P=P(w,v_0)$ embeds in $G$, the image $\bar{P}$ of $P$ is isomorphic to $P$, and the intersection of $\bar{P}$ and $F$ is empty. In addition, we can replace $P(w, v_0)$ by $C(w,v_0)$ if $w$ is weakly cyclically reduced and $w\notin \langle \beta \rangle$.

\end{cor}

\begin{proof}
The construction of the extension consists of choosing some finite points in $X$. Therefore, it is enough to choose all considering points far enough outside of $F$.
\end{proof}

\subsubsection{Property (FF)}\label{section-FF}

\begin{figure}[H]
\centering
\includegraphics[width=9.5cm]{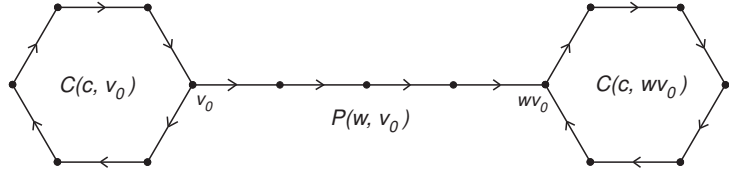}
\caption{The graph $Q_0=C(c, v_0)\cup P(w, v_0)\cup C(c, wv_0)$}
\label{model-Q_0}
\end{figure}

\begin{figure}
\centering \includegraphics[width=8cm]{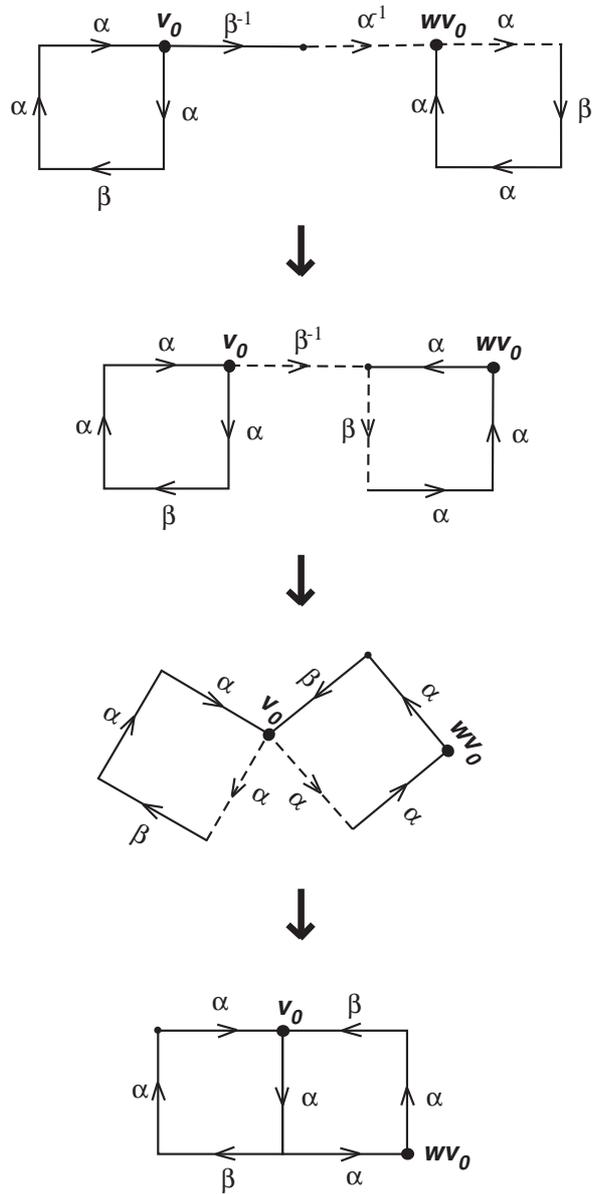}
\caption{Example of gluing double edges}
\label{gluing}
\end{figure}

Let $c=c_m\cdots c_1 $ be a weakly cyclically reduced word, such that $c\notin \langle\beta\rangle$. Let $w=w_k\cdots w_1$ be a reduced word, such that $w\notin \langle c \rangle$.
Let $C(c, v_0)$ be the cycle defined in Definition \ref{def-cycle}. Let $P(w, v_0)$ be the path defined in Definition \ref{def-path} such that every vertex of $P(w,v_0)$ (other than $v_0$) is distinct from every vertex in $C(c,v_0)$. Let $wv_0$ be the endpoint of $P(w,v_0)$. Let $C(c, wv_0)$ be the cycle with $i(c_1)=t(c_m)=wv_0$, such that every vertex of $C(c,wv_0)$ (other than $wv_0$) is distinct from every vertex in $P(w, v_0)\cup C(c, v_0)$ (see Figure \ref{model-Q_0}). Let us denote by $Q_0$ the union of $C(c, v_0)$, $P(w, v_0)$ and $C(c, wv_0)$. In general, this finite labeled graph $Q_0$ is not well-labeled. However, by identifying the successive edges with the same initial vertex and the same label, $Q_0$ becomes a well-labeled graph $Q$ (See Figure \ref{gluing} for an example of the process).

In the end of the process of identification of ``double edges", $Q$ has fewer edges than $Q_0$; however, the cycle $C(c,v_0)$ and $C(c, wv_0)$ are not modified, in the sense that the ``shapes'' of $C(c, v_0)$ and $C(c,wv_0)$ in $Q_0$ are the same as in $Q$. In other word, the quotient map $Q_0 \twoheadrightarrow Q$ restricted to $C(c,v_0)$ and to $C(c,wv_0)$ is injective (each one separately).

By construction, in each process, the graph has the following property:

\bigskip

\noindent \textbf{Property (FF)}
\begin{enumerate}
\item[(1)] the starting point of $C(c, v_0)$ is equal to its endpoint which is $v_0$;
\item[(2)] the starting point of $P(w, v_0)$ is different from its endpoint;
\item[(3)] the starting point of $C(c,wv_0)$ is equal to its endpoint which is $wv_0$.
\end{enumerate}

The acronym (FF) stands for ``Faithfulness for $w$ and fixed points of $c$''. Notice that $(2)$ comes from the fact that $w\notin \langle c \rangle$. When this process is finished, $Q$ will be one of the following four types (Figure \ref{fourQ}) of well-labeled graph satisfying the property (FF):

\begin{figure}[H]
\centering \includegraphics[width=12.5cm]{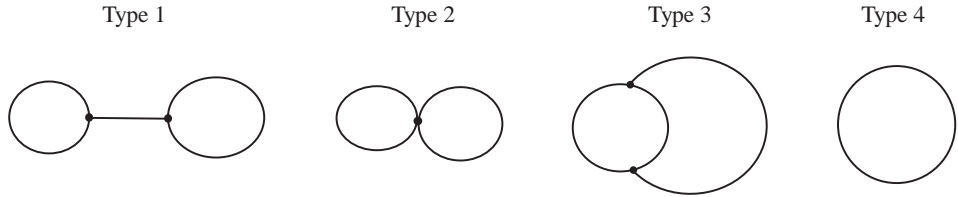}
\caption{Four types of $Q$}
\label{fourQ}
\end{figure}

\begin{prop}\label{extension}
For every one of the four types of well-labeled graph $Q=Q(c,w,v_0)$, there exists an extension $G$ of $G_0$ such that $Q$ embeds in $G$ and the image $Q(c,w,z_0)$ of $Q$ by the embedding has the property (FF), i.e. there exists $\alpha$ such that the word $w$ satisfies
$$
\left\{\begin{array}{l} cz_0=z_0; \\
wz_0 \neq z_0; \\
cwz_0=wz_0. \end{array} \right.
$$
where $z_0$ is the image of $v_0$ in $G$.
\end{prop}

We have to prove that every cycle in each types in Figure \ref{fourQ} contains at least one directed edge labeled by $\alpha$ or $\alpha^{-1}$. This is clear for the type 1, 2 and 4 since they have at most 2 cycles that represent $C(c,v_0)$ and $C(c,wv_0)$, and $c\notin \langle \beta \rangle$ by assumption. For the type 3, we can read around two subgraphs representing $C(c,v_0)$ and $P(w,v_0)^{-1}C(c,wv_0)P(w,v_0)$ from the vertex $v_0$. The labeling of the graph $P(w,v_0)^{-1}C(c, wv_0)P(w,v_0)$ is $w^{-1}cw$.

Let us recall the well-known theorem concerning the test for conjugacy of two words (see Theorem 1.3 in \cite{MaKaSol}):
\begin{thm}\label{conjugate}
Two words in the free group $\mathbb{F}_n$ define conjugate elements of $\mathbb{F}_n$ if and only if their cyclic reductions in $\mathbb{F}_n$ are cyclic permutations of one another.
\end{thm}

\begin{lem}\label{liste}
Let $c$ be a weakly cyclically reduced word, such that $c\notin \langle\beta\rangle$. Let $w$ be a reduced word, such that $w\notin \langle c \rangle$. If $c$ has the form $\gamma \beta^{l}$ with $\gamma \notin \langle \beta \rangle$, then $w^{-1}cw$ cannot be reduced to neither the form $\gamma \beta^{-k}$, nor the form $\gamma^{-1}\beta^k$ with $sign(k)=sign(l)$, $\forall k \in \mathbb{Z}$.
\end{lem}

\begin{proof}
Let $\gamma\beta^l$ with $\gamma= \gamma_n \cdots \gamma_1\notin \langle \beta \rangle$. By contradiction, let us suppose that $\gamma_n \cdots \gamma_1\beta^l$ is conjugate to $\gamma_n\cdots \gamma_1\beta^{-k}$ with $k,l>0$. Without loss of generality, we can suppose that $\gamma_1$, $\gamma_n \notin \{ \beta^{\pm 1} \}$. There are four types of cyclic permutations of $\gamma_n \cdots \gamma_1\beta^l$, which are $\gamma_n \cdots \gamma_1\beta^l$; $\beta^l\gamma_n \cdots \gamma_1$; $\beta^{l_1} \gamma_n \cdots \gamma_1 \beta^{l_2}$ with $l_1 + l_2 =l$; and $\gamma_p \cdots \gamma_1\beta^l \gamma_n \cdots \gamma_{p+1}$ for a certain $1\leq p \leq n$. Obviously, $\gamma_n\cdots\gamma_1\beta^{-k}$ cannot be of the first three types; so let us suppose that there exists $1\leq p\leq n$ such that $\gamma_p \cdots \gamma_1\beta^l \gamma_n \cdots \gamma_{p+1}=\gamma_n\cdots \gamma_1\beta^{-l}$ (since the two conjugate elements have the same length). By identification of the $l^{\textrm{th}}$ letter on the right of the two words, we have $\beta^{-1}=\gamma_{p+l}=\gamma_j$, for every $j$ multiple of $p+l$ modulo $n+l$, so in particular $\beta^{-1}=\gamma_{n-p}$. However, by identifying the $(n-p+l)^{\mathrm{th}}$ letter, which is $\beta$ for the left side, and $\gamma_{n-p}$ for the right side, we have $\beta=\gamma_{n-p}$ which contradicts with the first identification. The second case can be treated similarly.
\end{proof}

\noindent \textit{Proof of Proposition  \ref{extension}.} $\,$ As we mentioned before, it remains us to consider the type 3.

\begin{figure}[H]
\centering \includegraphics[width=4cm]{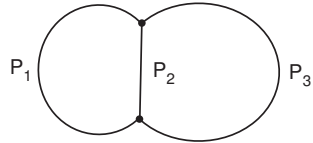}
\caption{Type 3 of $Q$}
\label{type3}
\end{figure}

In this graph, there are three cycles $C=P_1 \cup P_2$, $P_2 \cup P_3$ and $P_1\cup P_3$.

\smallskip
\noindent $\cdot$ \textbf{Claim.} If one of the three paths $P_1$, $P_2$ and $P_3$ has only edges labeled by $\beta^{\pm 1}$, then the other two paths both contains edges labeled by $\alpha^{\pm 1}$.

The claim allows to conclude. In fact, without loss of generality, suppose that $P_1$ has only edges labeled by $\beta^{\pm 1}$ and $P_2 \notin \langle \beta\rangle$ and $P_3\notin\langle\beta\rangle$. We first take an extension $G_1\supset G_0$ such that the image of $P_1$ is a path in $G_1$. Then we take an extension $G_2\supset G_1$ such that $P_2$ is a path in $G_2$ which connects the starting point and the endpoint of $P_1$ outside of the finite subset $P_1$; that is possible since the graph is well-labeled and $P_2$ contains edges labeled by $\alpha$. Finally, we take an extension $G_3 \supset G_2$ such that $P_3$ is a path in $G_3$ joining these two points outside of $P_1\cup P_2$.

We now prove the claim. Indeed, if two of these three paths were labeled by $\beta^{\pm 1}$, then $c$ would be the form of $\gamma \beta^l$ up to cyclic permutation and $w^{-1}cw$ would be the form of $\gamma \beta^{-k}$ or $\gamma^{-1}\beta^k$ with $sign(l)=sign(k)$ up to cyclic permutation, which contradicts with Lemma \ref{liste}.

\begin{flushright}
$\square$
\end{flushright}

\begin{cor}\label{cor-final}
Let $Q=Q(c,w,v_0)$ be a well-labeled graph. Let $F\subset G_0$ be a finite subset of $X$. There exists an extension $G$ of $G_0$ such that the image $Q(c,w,z_0)$ of $Q(c,w,v_0)$ in $G$ preserve the property (FF), and the intersection of $Q(c,w,z_0)$ and $F$ is empty.
\end{cor}

\begin{flushright}
$\square$
\end{flushright}

\subsubsection{Proof of Proposition  \ref{F_2-fidel1}}
Let $c=\alpha^{a_1}\beta^{b_1}\cdots \alpha^{a_n}\beta^{b_n}$ be a weakly cyclically reduced word on $\{\alpha^{\pm 1} $, $\beta^{\pm 1} \}$ (the other three types are similar). Let $w \in \langle\alpha, \beta \rangle \setminus \langle c \rangle$ be a reduced word on $\{\alpha^{\pm 1} $, $\beta^{\pm 1} \}$. We shall prove that the set
\begin{eqnarray}
\mathcal{V}_w=\{\alpha \in Sym(X) & \mid & \textrm{there exists a finite number of $x\in X$ such that} \nonumber \\
&& \textrm{$cx=x$, $cwx=wx$, and $wx \neq x$ }\} \nonumber
\end{eqnarray}
is meagre. For $K \subset X$ a finite subset of $X$, let
$$
V_{w,K}=\{\alpha \in Sym(X)\mid \big(\textrm{Fix($c$)}\cap w^{-1}\textrm{Fix($c$)}\cap \supp(w)\big)\subseteq K \},
$$
where $\supp(w)=\{x\in X \mid wx \neq x\}$.

The set $V_{w,K}$ is closed since if $\alpha_n$ converges to $\alpha$, then $c(\alpha_n,\beta)$ converges to $c(\alpha, \beta)$ and $w(\alpha_n, \beta)$ converges to $w(\alpha, \beta)$.
We shall prove that the interior of $V_{w,K}$ is empty.

\begin{lem}\label{def-on-F}
Let $\alpha'\in Sym(X)$ and $F \subset X$ be a finite subset of $X$. There exists $\alpha \in Sym(X)$ such that $\alpha |_{F}=\alpha' |_{F}$ and $\mathrm{supp}(\alpha)\subset F \cup \alpha'(F)$.
\end{lem}
\begin{proof}
Let us partition $F$ into finitely many pieces $F=\sqcup_{i=1}^m P_i$ according to the orbits of $\alpha'$. If $\alpha'(P_i)=P_i$, then define $\alpha |_{P_i}=\alpha'|_{P_i}$. If not, write $P_i=\{p_i$, $\alpha'(p_i)$, $\dots$, $\alpha'^{k_i}(p_i) \}$ with $\alpha'^{k_i+1}(p_i)\notin F$. Then define $\alpha |_{P_i}=\alpha'|_{P_i}$ and $\alpha(\alpha'^{k_i+1}(p_i))=p_i$.
\end{proof}

We see $X$ as the pre-graph $G_0$, where the $\beta^{\pm 1}$-edges of $G_0$ are seen as the transitive action of $\beta^{\pm 1}$ on $X$, which is fixed from the beginning.

Let $\alpha'\in V_{w,K}$ and let $F\subset X$ be a finite subset of $X$. Let $Y=F \cup \alpha'(F) \cup K$ be a finite subset of $X$. We construct a well-labeled graph $Q(c,w,v_0)$ as in Section \ref{section-FF}. We choose $z_0\notin Y$ and take $\alpha$ which is defined on $F$ as in Lemma \ref{def-on-F}, and which satisfies the property (FF) without touching any point of $Y$ (Corollary \ref{cor-final}). Consequently, $\alpha \notin V_{w,K}$ and $\alpha |_{F}=\alpha' |_{F}$.

\subsubsection{Proof of Proposition \ref{F_2-fidel2}}
We want to prove that for all $k\in \mathbb{Z}\setminus \{0\}$, the set
$$
\mathcal{V}_k =\{\alpha \in Sym(X) \mid c^k=\mathrm{Id}\}
$$
is closed and of empty interior.

Indeed, it is clearly closed. Moreover, let $\alpha' \in \mathcal{V}_k$ and let $F \subset X$ be a finite subset of $X$. Let $P(c^k,v_0)$ be the path defined in Definition \ref{def-path}. We choose $z_0\notin F\cup \alpha' (F)=:Y$ and take $\alpha$ which is defined on $F$ as in Lemma \ref{def-on-F}, and such that $P(c^k,z_0)$ is a path in $X$ not touching any point of $Y$. By consequent, $\alpha\notin \mathcal{V}_k$ and $\alpha |_{F}=\alpha' |_{F}$.

\subsection{Proof of Proposition \ref{F_2-moy}}
Let $c$ be a special word. Let $\{A_n\}_{n\geq 1}$ be a pairwise disjoint F\o lner sequence for $\beta$. Let $\{\varepsilon_l\}_{l\geq 1} > 0$ be a sequence tending to 0. Let us write
\begin{eqnarray}
\mathcal{U}_3=\bigcap_l \bigcap_{N\in\mathbb{N}} \{\alpha\in Sym(X)&\mid & \textrm{there exists $k\geq N$ such that $A_k\subset \mathrm{Fix}(c)$ and} \nonumber \\ && \textrm{$|A_k\vartriangle \alpha A_k|<\varepsilon_l |A_k|$ }\}.\nonumber
\end{eqnarray}
Set $\varepsilon_l = \varepsilon$. We want to prove that the set
$$
\mathcal{V}_N :=\{\alpha\in Sym(X) \mid \forall k\geq N, \textrm{ $A_k \nsubseteq \mathrm{Fix}(c)$ or $|A_k\vartriangle \alpha A_k|\geq\varepsilon |A_k|$ }\}
$$
is closed and of empty interior. We treat the case $c=\alpha^{a_1}\beta^{b_1}\cdots \alpha^{a_n}\beta^{b_n}$ (the other three types are similar). Let $M=\textrm{max}_j |b_j|$ and set
$$
E_k:=\cup_{i=-M}^M \beta^i(A_k),
$$
a finite set of $X$.

\bigskip

\noindent $\cdot$ $\mathcal{V}_N$ \textit{is closed.} $\quad$ Since $\mathcal{V}_N= \cap_{k\geq N}\mathcal{V}_{N,k}$ where
$$
\mathcal{V}_{N,k}:=\{\alpha\in Sym(X) \mid \textrm{ $A_k \nsubseteq \mathrm{Fix}(c)$ or $|A_k\vartriangle \alpha A_k|\geq\varepsilon |A_k|$ }\},
$$
it is enough to prove that $\mathcal{V}_{N,k}$ is closed. So let $\{\alpha_n\}_{n\geq 1}$ be a sequence in $\mathcal{V}_{N,k}$ which converges to $\alpha \in Sym(X)$. Since $E_k$ is finite, there exists $n_0$ such that $\alpha |_{E_k}=\alpha_n |_{E_k}$, $\forall n\geq n_0$. Therefore, $\alpha \in \mathcal{V}_{N,k}$ because $A_k \subset E_k$.

\bigskip

\noindent $\cdot$ $\mathcal{V}_N$ \textit{is of empty interior.} $\quad$ Let us distinguish two cases:

First, suppose that $S(\alpha)=S(\beta)=0$. Let $\alpha'\in \mathcal{V}_N$. Let $F \subset X$ be a finite subset of $X$. We choose $m \gg N$ such that $(F\cup \alpha'(F))\cap E_m=\emptyset$. We define $\alpha |_{E_m}=$Id and $\alpha |_F=\alpha' |_F$. Then $A_m \subset \mathrm{Fix}(c)$ since $S(\beta)=0$, and $|A_m\vartriangle \alpha A_m|=0$ since $\alpha(A_m)=A_m$. So $\alpha\notin \mathcal{V}_N$.

Second, suppose that $S(\alpha)$ divides $S(\beta)$. Let $\alpha'\in \mathcal{V}_N$. Let $F \subset X$ be a finite subset of $X$. We choose $m \gg N$ such that $(F\cup \alpha'(F))\cap E_m=\emptyset$ and $|A_m \vartriangle \beta^{-\frac{S(\beta)}{S(\alpha)}}(A_m)|<\varepsilon |A_m|$; this is possible as $\{A_m\}$ is a F\o lner sequence for $\beta$. We define
$$
\alpha(x)=\beta^{-\frac{S(\beta)}{S(\alpha)}}(x), \quad \forall x\in E_m,
$$
and $\alpha |_F=\alpha' |_F$. Then,
$$
c(x)=\beta^{-\frac{S(\beta)}{S(\alpha)}a_1}\beta^{b_1}\cdots \beta^{-\frac{S(\beta)}{S(\alpha)}a_n}\beta^{b_n}(x)=\beta^{-\frac{S(\beta)}{S(\alpha)}S(\alpha)}\beta^{S(\beta)}(x)=x,
$$
for all $x\in E_m$. In particular, $A_m \subset \mathrm{Fix}(c)$. In addition,
$$
|A_m \vartriangle \alpha A_m|=|A_m \vartriangle \beta^{-\frac{S(\beta)}{S(\alpha)}}(A_m)|<\varepsilon |A_m|,
$$
so $\alpha \notin \mathcal{V}_N$.

\subsection{Proof of Proposition \ref{F_2-transitive}}

The proof follows from the three claims:

\bigskip

\noindent $\cdot$ \textbf{Claim 1.} Let $G$ be a group and $H<G$ be a finite index subgroup of $G$. Then, for all $g\in G$, there exists $n\geq 1$ such that $g^n\in H$.

Indeed, let $N$ be the core of $H$, that is $N=\bigcap_{x\in G}x^{-1}Hx \subset H$. The subgroup $N$ is a finite index normal subgroup of $G$. Then for all $g\in G$, $g^m \in N$, where $m=[G:N]$.
\bigskip

\noindent $\cdot$ \textbf{Claim 2.} The set
$$
\mathcal{U}_5=\{\alpha\in Sym(X) \mid \textrm{$\forall n$, $m\in \mathbb{Z}\setminus \{0\}$, the $\langle\alpha^n,\beta^m\rangle$-action on $X$ is transitive } \}
$$
is in $\mathcal{U}_4$.

Indeed, let $\alpha\in \mathcal{U}_5$. Let $H < \langle\alpha,\beta\rangle$ be a finite index subgroup. Then by Claim 1, there exist $n_0$, $m_0$ such that $\alpha^{n_0}$ and $\beta^{m_0}$ are in $H$, so $\langle\alpha^{n_0},\beta^{m_0}\rangle < H$. Since the $\langle\alpha^{n_0},\beta^{m_0}\rangle$-action on $X$ is transitive by hypothesis, the $H$-action on $X$ is also transitive.

\bigskip

\noindent $\cdot$ \textbf{Claim 3.} The set $\mathcal{U}_4$ is generic in $Sym(X)$.

It is enough to prove that the set $\mathcal{U}_5$ is generic since $\mathcal{U}_5 \subset \mathcal{U}_4$. So let us prove that for all $n$ and $m$, the set $\mathcal{V}_{n,m}=\{\alpha\in Sym(X) \mid \langle\alpha^{n},\beta^{m}\rangle$-action on $X$ is not transitive $\}$ is closed and it has empty interior.

\bigskip

\noindent \textit{$\cdot$ $\mathcal{V}_{n,m}$ is closed.}
\begin{eqnarray}
\mathcal{V}_{n,m}&=& \{\alpha \in Sym(X) \mid \textrm{$\exists x$, $y \in X$ such that $\forall w \in \langle \alpha^n, \beta^m \rangle$, $wx\neq y$} \} \nonumber \\
&=& \{\alpha \in Sym(X) \mid \textrm{$\exists (x_i, x_j) \in S\times S$ such that $\forall w \in \langle \alpha^n, \beta^m \rangle$, $wx_i\neq x_j$}\} \nonumber \\
&=& \bigcup_{(x_i, x_j) \in S\times S} \{\alpha \in Sym(X) \mid \forall w \in \langle \alpha^n, \beta^m \rangle, \, wx_i\neq x_j \} \nonumber
\end{eqnarray}
where $S$ is a finite family of representatives for $\beta^m$-orbits. It is clear that the set $\{\alpha \in Sym(X) \mid \forall w \in \langle \alpha^n, \beta^m \rangle, \, wx_i\neq x_j \}$ is closed. So $\mathcal{V}_{n,m}$ is closed as a finite union of closed sets.

\bigskip

\noindent \textit{$\cdot$ $\mathcal{V}_{n,m}$ is of empty interior.} $\quad$ Let $\alpha'\in \mathcal{V}_{n,m}$ and let $F\subset X$ be a finite subset of $X$. Let $Y:=F\cup \alpha'(F)$ be a finite subset of $X$. We choose representatives for $\beta^m$-orbits outside of $Y$, and form a finite family $S=\{x_1$, $\dots$, $x_m\}$ of $X$; this is possible since the $\beta^m$-orbits are infinite. We define $\alpha$ on $F$ as in Lemma \ref{def-on-F}. Inductively on $1\leq i \leq m-1$, in each $\beta^m$-orbit $O(x_i)$ of $x_i$, we choose $n-1$ points $\{p_{i,1}$, $p_{i,2}$, $\dots$, $p_{i,n-1}\}$ outside of $Y$ and define
\begin{enumerate}
\item[$\cdot$] $\alpha(x_i)=p_{i,1}$;
\item[$\cdot$] $\alpha(p_{i,j})=p_{i, j+1}$, $ \forall 1 \leq j\leq n-2$;
\item[$\cdot$] $\alpha(p_{i,n-1})=x_{i+1}$.
\end{enumerate}
Then, in $O(x_m)$, we choose $n-1$ points $\{p_{m,1}$, $\dots$, $p_{m,n-1}\}$ outside of $Y$ and define
\begin{enumerate}
\item[$\cdot$] $\alpha(x_m)=p_{m,1}$;
\item[$\cdot$] $\alpha(p_{m,j})=p_{m, j+1}$, $ \forall 1 \leq j\leq n-2$;
\item[$\cdot$] $\alpha(p_{m,n-1})=x_1$.
\end{enumerate}
By construction, $\alpha^n(x_i)=x_{i+1}$, $\forall 1\leq i\leq m-1$, and $\alpha^n(x_m)=x_1$, so the $\langle \alpha^n, \beta^m \rangle$-action is transitive.

\section{Construction of $\mathbb{F}_2 \ast_\mathbb{Z} \mathbb{F}_2$}

Let $X$ be a countable infinite set. Let $c=c(\alpha,\beta)$ be a special word. Let $G:=\mathbb{F}_2=\langle \alpha$, $\beta \rangle$ be constructed as in Chapter 3. Let $\{A_n\}_{n=1}^{\infty}$ be a F\o lner sequence such that $c(A_n)=A_n$, $\forall n\geq 1$. Let $Z_c=\{\sigma \in Sym(X)\mid \sigma c = c \sigma\}$ be the centralizer of $c$. Let $\alpha'=\sigma^{-1}\alpha\sigma$, $\beta'=\sigma^{-1}\beta \sigma$, and let $H:=\langle \alpha'$, $\beta' \rangle$. Let $A=\langle c \rangle$ be the subgroup of $G$ generated by $c$. We consider $\mathbb{F}_2 \ast_\mathbb{Z} \mathbb{F}_2=G\ast_{A}H$ the amalgamated free product of $G$ and $H$ along $A$. For all $\sigma \in Z$, the action of $G\ast_A H$ on $X$ is given by $g\cdot x=g(\alpha,\beta)x=gx$, and $h\cdot x=h(\alpha', \beta')x=\sigma^{-1} h(\alpha,\beta)\sigma x=\sigma^{-1}h\sigma x$, for all $g\in G$ and $h\in H$.
\begin{lem}
The set $Z_c$ is closed in $Sym(X)$. In particular, $Z_c$ is a Baire space.
\end{lem}
\begin{proof}\label{amalgame-fidel}
The application $p:Sym(X)\rightarrow Sym(X);$ $\sigma \mapsto [\sigma, c]$ is continuous. So $Z_c=p^{-1}\{\textrm{Id} \}$ is closed since $\{\textrm{Id}\}$ is closed in $Sym(X)$.
\end{proof}
\begin{prop}
The set
$$
\mathcal{O}_1=\{\sigma\in Z_c \mid \textrm{ the action of $G\ast_A$H on $X$ is faithful } \}
$$
is generic in $Z_c$.
\end{prop}

\begin{proof}
For all $w\in G\ast_A H$, let us denote by $w^{\sigma}$ the corresponding element of $Sym(X)$ given by the above action, i.e. if $w=ag_n h_n\cdots g_1 h_1$, with $a\in A$, $g_i\neq e\in G\setminus A$ and $h_i \neq e\in H \setminus A$, for all $i$, then
 $$w^{\sigma}=ag_n\sigma^{-1}h_n\sigma\cdots g_1 \sigma^{-1}h_1\sigma.$$
We want to prove that the set
$$
\mathcal{O}_1=\bigcap_{w\neq e \in G\ast_A H} \{\sigma\in Z_c \mid \exists x\in X \textrm{ such that $w^{\sigma}x \neq x$ } \}
$$
is generic in $Z_c$. Therefore, we shall prove that the set
$$
\mathcal{V}_w=\{\sigma \in Z_c \mid w^{\sigma}=\mathrm{Id}_X\}
$$
is closed and of empty interior in $Z_c$.

The set $\mathcal{V}_w$ is closed in $Z_c$ because the application $Z_\sigma \rightarrow Sym(X)$; $\sigma \mapsto w^{\sigma}$ is continuous.

To see that the set $\mathcal{V}_w$ is of empty interior, let $\sigma'\in \mathcal{V}_w$, and let $F\subset X$ be a finite subset of $X$. Notice that if $F=F_1 \sqcup F_2$ with $F_1\subset $ Fix($c$) and $F_2\cap \Fix(c)=\emptyset$, then $\sigma'(F_1) \subset $ Fix($c$) and $\sigma'(F_2)\cap \Fix(c)=\emptyset$ because $\sigma'(\textrm{Fix$(c)$})=$ Fix$(c)$, for all $\sigma'\in Z_c$. So we define $\sigma |_{F_1}=\sigma' |_{F_1}$ as in Lemma \ref{def-on-F}, and $\sigma |_{X\setminus \mathrm{Fix}(c)} = \sigma' |_{X\setminus \mathrm{Fix}(c)}$. Therefore, we have defined $\sigma$ on $Y:=(F\cup \sigma'(F))\cup (X\setminus \mathrm{Fix}(c))$, and $\sigma |_Y $ commutes with $c |_Y$. Let us now define $\sigma $ on $X\setminus Y$ in a way that $\sigma \in Z_c \setminus \mathcal{V}_w$. For all $g\in G \setminus A$ and $h\in H\setminus A$, let
$$
\widehat{g}=\{x\in X \mid cx=x, \, cgx=gx \textrm{ and $gx\neq x$ }\},
$$
$$
\widehat{h}=\{x\in X \mid cx=x, \, chx=hx \textrm{ and $hx\neq x$ }\}.
$$
Recall that we are considering the word $w^{\sigma}=ag_n\sigma^{-1}h_n\sigma\cdots g_1 \sigma^{-1}h_1\sigma$. Choose any $x_0 \in X\setminus Y$. By induction on $1\leq i\leq n$, we choose $x_{4i-3}\in \widehat{h_i}$ such that $x_{4i-3}$ is different from the finite set of points $x_1$, $\dots$, $x_{4i-4}$ chosen until the $(i-1)^{\mathrm{th}}$ step. This is possible since $\widehat{h_i}$ is infinite by Proposition \ref{F_2-fidel1}. Then we define $\sigma x_{4i-4}:=x_{4i-3}$ and $\sigma x_{4i-3}:=x_{4i-4}$. This is well-defined because $x_{4i-4}$, $x_{4i-3}\in \mathrm{Fix}(c)$. We set $h_i x_{4i-3}=: x_{4i-2}$ which is different from $x_{4i-3}$ and which is fixed by $c$, by definition of $\widehat{h_i}$. We choose $x_{4i-1}\in \widehat{g_i}$ such that $x_{4i-1}$ is different from the finite set of points chosen so far. This is again possible since $\widehat{g_i}$ is infinite (Proposition \ref{F_2-fidel1}). Then we define $\sigma x_{4i-2}:=x_{4i-1}$ and $\sigma x_{4i-1}:=x_{4i-2}$. This is also well-defined because $x_{4i-2}$, $x_{4i-1}\in \mathrm{Fix}(c)$. We finally set $g_i x_{4i-1}=:x_{4i}$. By construction, the $4n$ points defined by the subwords on the right of $w^{\sigma}$ are all distinct. In particular, $w^{\sigma}x_0=ax_{4n}=x_{4n}\neq x_0$. Besides, this construction works also for the other three types of word $w$ since we are treating all subwords of $w$. At last, if $w=g\in G\setminus \{\mathrm{Id}\}$, then there exists $x\in X$ such that $gx\neq x$ since $G$ acts faithfully on $X$. Therefore, $\sigma$ constructed in this way is beautifully in $Z_c\setminus \mathcal{V}_w$ and $\sigma' |_F=\sigma |_F$.

\end{proof}

\begin{prop}\label{amalgame-moy}
The set
\begin{eqnarray}
\mathcal{O}_2=\{\sigma \in Z_c &\mid &\textrm{ there exists $\{A_{n_k}\}_{k\geq 1}$ a subsequence of $\{A_n\}_{n\geq 1}$ such that }\nonumber \\
 && \textrm{ $\sigma(A_{n_k})=A_{n_k}$, $\forall k\geq 1$ }\} \nonumber
\end{eqnarray}
is generic in $Z_c$.
\end{prop}

\begin{proof}
We want to prove that the set
$$
\mathcal{O}_2=\bigcap_{N\in \mathbb{N}}\{\sigma\in Z_c \mid \exists n\geq N \textrm{ such that $\sigma(A_n)=A_n$ }\}
$$
is generic in $Z_c$. So we shall prove that the set
$$
\mathcal{V}_N=\{\sigma \in Z_c \mid \forall n\geq N, \, \sigma(A_n)\neq A_n \}
$$
is closed and of empty interior in $Z_c$.

\bigskip

\noindent \textit{$\cdot$ $\mathcal{V}_N$ is closed.} $\quad$ It is enough to prove that the set
$$
V_{n,N}=\{\sigma\in Z_c \mid \sigma(A_n)\neq A_n\}
$$
is closed since $\mathcal{V}_N=\cap_{n\geq N}V_{n,N}$. Let $\{\sigma_m\}_{m\geq 1}\subset V_{n,N}$ be a sequence converging to $\sigma$ in $Z_c$. Since $A_n$ is finite, there exists $m_0$ such that $\sigma_m(A_n)=\sigma(A_n)$, $\forall m\geq m_0$. Thus we have $\sigma(A_n)\neq A_n$ since $\sigma_m(A_n)\neq A_n$.

\bigskip

\noindent \textit{$\cdot$ $\mathcal{V}_N$ is of empty interior.} $\quad$ Let $\sigma'\in \mathcal{V}_N$ and let $F\subset X$ be a finite subset of $X$. Let $Y:=(F\cup \sigma'(F))\cup(X\setminus \mathrm{Fix}(c))$. Since $A_n \subset \mathrm{Fix}(c)$ (Proposition \ref{F_2-moy}), there exists $n\geq N$ such that $A_n\cap Y=\emptyset$. We take then $\sigma\in Z_c$ which fixes $A_n$ and $\sigma |_Y=\sigma' |_Y$. Therefore, $\sigma \in Z_c \setminus \mathcal{V}_N$ and $\sigma |_F=\sigma' |_F$.
\end{proof}

Let $\sigma \in \mathcal{O}_1\cap \mathcal{O}_2$. Let $\{A_{n_k}\}_{k\geq 1}$ be a subsequence of $\{A_n\}_{n\geq 1}$ such that $\sigma(A_{n_k})=A_{n_k}$, $\forall k\geq 1$. We claim that $\{A_{n_k}\}_{k\geq 1}$ is a F\o lner sequence for $G\ast_A H$. Indeed, for all $g\in G$ and for all $h\in H$, we have
$$
\lim_{k\rightarrow \infty}\frac{|A_{n_k} \vartriangle g\cdot A_{n_k}|}{|A_{n_k}|} = \lim_{k\rightarrow \infty}\frac{|A_{n_k} \vartriangle g(\alpha,\beta)A_{n_k}|}{|A_{n_k}|}=0,
$$
\begin{eqnarray}
\lim_{k\rightarrow \infty}\frac{|A_{n_k} \vartriangle h\cdot A_{n_k}|}{|A_{n_k}|}&=&\lim_{k\rightarrow \infty}\frac{|A_{n_k} \vartriangle h(\alpha',\beta') A_{n_k}|}{|A_{n_k}|}\nonumber \\
&=&\lim_{k\rightarrow \infty}\frac{|A_{n_k} \vartriangle \sigma^{-1}h(\alpha,\beta)\sigma A_{n_k}|}{|A_{n_k}|}\nonumber \\
&=& \lim_{k\rightarrow \infty}\frac{|\sigma A_{n_k} \vartriangle h(\alpha,\beta)\sigma A_{n_k}|}{|A_{n_k}|}\nonumber \\
&=&\lim_{k\rightarrow \infty}\frac{|A_{n_k} \vartriangle h(\alpha,\beta)A_{n_k}|}{|A_{n_k}|}=0,\nonumber
\end{eqnarray}
since $\{A_{n_k}\}$ is F\o lner for $G$ and $\sigma(A_{n_k})=A_{n_k}$. Therefore, we have:

\begin{thm}\label{amalgam}
There exists a transitive, faithful and amenable action of the group $\langle \alpha,\beta\rangle \ast_{\langle c \rangle} \langle \alpha',\beta'\rangle$ on $X$.
\end{thm}

\begin{lem}
Let $c=c(\alpha,\beta)$ be any word (not necessarily special) on $\{\alpha^{\pm 1},\beta^{\pm 1}\}$. There exists an automorphism $a$ of $\mathbb{F}_2$ such that $a(c)$ is a special word.
\end{lem}

\begin{proof}
Let us recall some properties of automorphisms of free groups. The reader can find more details in \cite{LynSch}. Let $\mathbb{F}_n$ be a free group with a finite basis $X$ of $n$ elements. We consider the following endomorphisms of $\mathbb{F}_n$. For any $x\in X$, let $\varphi_x$ be the endomorphism defined by $\varphi_x:$ $x\mapsto x^{-1}$; $y\mapsto y$, $\forall y\in X\setminus \{x\}$. For any $x\neq y\in X$, let $\psi_{xy}:$ $x\mapsto xy$; $z\mapsto z$, $\forall z\in X\setminus \{x\}$. In both cases, the image of $X$ is another basis for $\mathbb{F}_n$, and $\varphi_x$ and $\psi_{xy}$ are automorphisms of $\mathbb{F}_n$, called the Nielsen generators for Aut$(\mathbb{F}_n)$, and they generate Aut$(\mathbb{F}_n)$. Let $\mathbb{F}_n/\mathbb{F}_n'\simeq \mathbb{Z}^n$ be the abelianization of $\mathbb{F}_n$. We have Aut$(\mathbb{Z}^n)\simeq GL_n(\mathbb{Z})$. The Nielsen generators for Aut$(\mathbb{F}_n)$ induce the following generators for Aut($\mathbb{Z}^n$):
$$
\overline{\varphi}_x:x\mapsto -x; \quad y\mapsto y, \quad \forall y \in X\setminus\{x\};
$$
$$
\overline{\psi}_{xy}:x\mapsto x+y; \quad z\mapsto z, \quad \forall z \in X\setminus\{x\}.
$$
Thus, we conclude that the natural maps from Aut$(\mathbb{F}_n)$ into Aut$(\mathbb{Z}^n)$ is an epimorphism. Notice that for a word $c$ to be a special word depends only on its image in $\mathbb{Z}^2$. Therefore, in order to prove the Lemma, it is enough to find a matrix $M\in GL_2(\mathbb{Z})$ such that the exponent sum $S(\alpha)':=S_{a(c)}(\alpha)$ of exponents of $\alpha$ in the word $a(c)$ divides the exponent sum $S(\beta)':=S_{a(c)}(\beta)$ of exponents of $\beta$ in the word $a(c)$, where $a\in \mathrm{Aut}(\mathbb{F}_2)$ is a reciprocal image of $M$ by the epimorphism Aut$(\mathbb{F}_2)\rightarrow \mathrm{Aut}(\mathbb{Z}^2) $. In fact, once we have $c=c(\alpha,\beta)$ with $S(\alpha)$ dividing $S(\beta)$, we can obtain a weakly cyclically reduced word by conjugating $c$, and the conjugation is an automorphism of $\mathbb{F}_2$.

If $S(\beta)=0$, $c$ is already a special word. If $S(\alpha)=0$ and $S(\beta)\neq 0$, then we apply the matrix $\left( \begin{array}{cc} 0 & 1\\ 1 & 0 \end{array} \right) \in GL_2(\mathbb{Z})$ which exchanges $S(\alpha)$ and $S(\beta)$. So suppose that $S(\alpha)\neq 0 \neq S(\beta)$. Let $d=\mathrm{gcd}(S(\alpha),S(\beta))$ be the greatest common divisor of $S(\alpha)$ and $S(\beta)$. By Bézout's identity, there exist relatively prime integers $p$, $q$ such that $pS(\alpha)+qS(\beta)=d$. Since $\mathrm{gcd}(p,-q)=1$, there exist $r$, $t$ such that $rp-tq=1$ again by Bézout's identity. Then, the matrix $M=\left( \begin{array}{cc} p & q\\ t & r \end{array} \right)$ is in $GL_2(\mathbb{Z})$ and it sends $\left( \begin{array}{c} S(\alpha) \\ S(\beta) \end{array} \right)$ to $\left( \begin{array}{c} d \\ tS(\alpha)+rS(\beta) \end{array}\right)$. Therefore, $S(\alpha)'=d$ divides $S(\beta)'=tS(\alpha)+rS(\beta)$.
 \end{proof}
From Theorem \ref{amalgam} and the previous Lemma, we have:

\begin{thm}\label{amalgam_any_c}
Let $c=c(\alpha,\beta)$ be any word on $\{\alpha^{\pm 1},\beta^{\pm 1}\}$. Then the group $\langle \alpha,\beta\rangle \ast_{\langle c \rangle} \langle \alpha',\beta'\rangle$ admits a transitive, faithful and amenable action.
\end{thm}

A result of G. Baumslag \cite{Baum} shows that these groups are residually finite.

\smallskip

Furthermore, let $H$ be a finite index subgroup of $\mathbb{F}_2 \ast_{\mathbb{Z}} \mathbb{F}_2$. Then $K:=H\cap \mathbb{F}_2$ is a finite index subgroup of $\mathbb{F}_2$ so that the $H$-action on $X$ is transitive since the $K$-action is transitive by Proposition \ref{F_2-transitive}. Therefore, we have:

\begin{thm}\label{amalgam-index-in-A}
For any finite index subgroup $H$ of $\langle \alpha,\beta\rangle \ast_{\langle c \rangle} \langle \alpha',\beta'\rangle$, $H$ admits a transitive, faithful and amenable action.
\end{thm}

\section{Applications}

Let us recall the class of all countable groups appeared in \cite{GlMo}:
$$
\mathcal{A}=\{\textrm{ $G$ countable $\mid$ $G$ admits a faithful transitive amenable action }\}.
$$

Let $\Sigma_g$ be a closed oriented surface of genus $g\geq 2$. It is well-known that the fundamental groups $\Gamma_g=\pi_1 (\Sigma_g)$ of $\Sigma_g$ has a presentation
$$
\pi_1(\Sigma_g)=\langle \textrm{$a_1$, $b_1$, $\dots$, $a_g$, $b_g$ } \mid \prod_{i=1}^g [a_i,b_i]\, \rangle.
$$
In particular, we have $\pi_1(\Sigma_2)= \langle a_1,b_1\rangle \ast_{\langle c \rangle} \langle a_2,b_2\rangle$ where $c=[a_1,b_1]=[a_2,b_2]$. Therefore, $\pi_1(\Sigma_2)\in \mathcal{A}$ by Theorem \ref{amalgam_any_c} (or already by Theorem \ref{amalgam} since $c=[a_1,b_1]$ is a special word).
Now, let $\Sigma_g$ be a closed oriented surface of genus $g\geq 3$. Viewing $\Sigma_g$ as $(g-1)$ tori glued on a central one, the cyclic group $\mathbb{Z}/(g-1)\mathbb{Z}$ acts properly and freely on $\Sigma_g$, and the quotient space is $\Sigma_2$. Therefore $\pi_1(\Sigma_g)$ injects into $\pi_1(\Sigma_2)$ as a subgroup of index $(g-1)$ (in other words, $\Sigma_g$ is a $(g-1)$-sheeted regular covering of $\Sigma_2$). Consequently, $\pi_1(\Sigma_g)$ is in $\mathcal{A}$ by Theorem \ref{amalgam-index-in-A}. Moreover, the fundamental group of a torus $\pi_1(\mathbb{T}^2)=\pi_1(\Sigma_1)$ is isomorphic to $\mathbb{Z}^2$, an amenable group. Therefore, we have:

\begin{thm}\label{Sigma_g}
Let $\Sigma_g$ be a closed oriented surface of genus $g\geq 1$. The fundamental group $\Gamma_g=\pi_1 (\Sigma_g)$ of $\Sigma_g$ admits a transitive, faithful and amenable action, for all $g\geq 1$.
\end{thm}

\begin{cor}
For any compact surface $S$, the fundamental group $\pi_1 (S)$ is in $\mathcal{A}$.
\end{cor}

\begin{proof}
First of all, we can suppose that $S$ is oriented. In fact, it is well-known that if $S$ is a non-oriented connected surface, then there exists a oriented 2-sheeted covering space $\tilde{S}$ (cf. \cite{Ful}). Then $\pi_1(\tilde{S})$ is a subgroup of index 2 of $\pi_1(S)$ so that it is \textit{co-amenable} in $\pi_1(S)$ (a subgroup $H<G$ is co-amenable if the $G$-action on $G/H$ is amenable). Therefore, in order that $\pi_1(S)\in \mathcal{A}$, it suffices to have $\pi_1(\tilde{S})\in \mathcal{A}$ by Proposition 1. (vi) in \cite{GlMo}.

 If $S$ is a closed oriented surface (i.e. without boundary), $S$ is either a sphere or a finite connected sum of tori $\Sigma_g$, $g\geq 1$; so $\pi_1(S)\in \mathcal{A}$ in both cases.
 If $S$ is a surface with boundary components, then $\pi_1(S)$ is a free group (the fundamental group of a sphere with $p$ boundary components is a free group of rank $p-1$, and the fundamental group of $\Sigma_g$ with $p$ boundary components is a free group of rank $2g+p-1$, $\forall g\geq 1$), so it is again in $\mathcal{A}$ by van Douwen's theorem.

\end{proof}

\begin{ex}\label{bundle}
\textbf{ Surface bundles over $\mathbb{S}^1$}

A surface bundle over $\mathbb{S}^1$ is a closed 3-manifold which is constructed as a fiber bundle over the circle with fiber a closed surface. The fundamental group $G$ of such bundle can be viewed as an HNN-extension
$$
G=\pi_1(M_{\phi})=\langle \Gamma_g, t \mid txt^{-1}=\phi_{\ast}(x), \,\forall x\in \Gamma_g\rangle,
$$
where $\phi: \Sigma_g \rightarrow \Sigma_g$ is a homeomorphism. Thus, we have a short exact sequence:
$$
0 \rightarrow \Gamma_g \rightarrow G \rightarrow \mathbb{Z} \rightarrow 0.
$$
The subgroup $\Gamma_g$ is co-amenable in $G$ since it is normal in $G$ and $G/\Gamma_g \simeq \mathbb{Z}$ is amenable. Therefore, we have $G \in \mathcal{A}$.

The Thurston's virtual fibration conjecture states that \cite{Thr}:

\begin{quote}
\emph{Every closed, irreducible, atoroidal 3-manifold $M$ has a finite-sheeted cover which fibres over the circle.}
\end{quote}

It follows from the conjecture that the fundamental group $\pi_1(M)$ is in $\mathcal{A}$ since it contains a finite index subgroup which is in $\mathcal{A}$.
\end{ex}

\bibliographystyle{amsplain}
\bibliography{amalgambib}

\smallskip

\begin{quote}
Soyoung \textsc{Moon} \\
Institut de Math\'ematiques \\
Universit\'e de Neuch\^atel\\
11, Rue Emile Argand - BP 158\\
2009 Neuch\^atel - Switzerland

\smallskip

E-mail: \url{so.moon@unine.ch}
\end{quote}
\end{document}